# Fixed point theorems and convergence theorems for a generalized nonexpansive mapping in uniformly convex Banach spaces.

Chang Il Rim, Jong Gyong Kim

# Fixed point theorems and convergence theorems for a generalized non-expansive mapping in uniformly convex Banach spaces.

Chang Il Rim, Jong Gyong Kim

Faculty of Mathematics, Kim Il Sung University, DPR Korea

Faculty of Mathematics, Kim Il Sung University, DPR Korea

**Abstract** In this paper, we prove the existence of fixed points of mappings satisfying the condition ($D_a$), a kind of generalized nonexpansive mappings, on a weakly compact convex subset in a Banach space satisfying Opial's condition. And we use Sahu([6]) and Thakur([10])'s iterative scheme to establish several convergence theorems in uniformly convex Banach spaces and give an example to show that this scheme converges faster than the scheme in [1].

## 1. Introduction

In many literatures, several iterative schemes for fixed points of mappings have been studied.

In 2000, Noor([4]) introduced the following iterative scheme for general variational inequalities and studied the convergence criteria of this scheme.

$$\begin{cases} x = x_1 \in C \\ x_{n+1} = (1-a_n)x_n + a_n T y_n \\ y_n = (1-b_n)x_n + b_n T z_n \\ z_n = (1-c_n)x_n + c_n T x_n \end{cases}, n \in \mathbf{N} \quad (*)$$

, where $\{a_n\}, \{b_n\}, \{c_n\} \subset (0,1)$.

In 2016, Sahu([6]) and Thakur([10]) introduced the following iterative scheme for nonexpansive mappings in uniformly convex Banach spaces and claimed that this scheme converges to a fixed point of a contraction mapping faster than all the known iterative schemes.

$$\begin{cases} x = x_1 \in C \\ x_{n+1} = (1-a_n)Tz_n + a_n T y_n \\ y_n = (1-b_n)z_n + b_n T z_n \\ z_n = (1-c_n)x_n + c_n T x_n \end{cases}, n \in \mathbf{N} \quad (**)$$

, where $\{a_n\}, \{b_n\}, \{c_n\} \subset (0,1)$.

On the other hand, in recent years some new generalized nonexpansive mappings have been introduced and fixed point theorems and convergence theorems for these mappings have been studied.

In 2008, Suzuki([9]) defined a new kind of generalized nonexpansive mappings, ($C$)-mapping, and a lot of existence theorems and convergence theorems for a fixed point of ($C$)-mappings have been developed.([2]) In 2019, [3] introduced (**) for ($C$)-mappings in uniformly convex Banach spaces to prove the convergence theorems to a fixed point and gave a numerical example to show that this scheme converges faster than other schemes.

In 2018, [1] intoduced another new kind of generalized nonexpansive mappings, ($D_a$)-mappings and applied (*) to these mappings with $1/2 < a \leq a_n = b_n = c_n \leq b < 1$ to establish convergence theorems in Banach spaces.

In this paper we establish the existence of a fixed point of ($D_a$)-mappings in weakly compact convex subsets of Banach spaces and the convergence theorems of (**) for a fixed point of these mappings. We also give an example to show that (**) converges faster than (*).

Here are some useful concepts and properties.

**Definition 1** ([5]) Let $T: C \to C$ be a mapping on a subset $C$ of a Banach space $X$. If $F(T) \neq \phi$ and $T$ satisfies the following inequality for all $x \in C, p \in F(T)$, then $T$ is said to be quasinonexpansive, where $F(T)$ is the set of fixed points of $T$.

$$\|Tx - p\| \leq \|x - p\|.$$

**Definition 2** ([9]) Let $T: C \to C$ be a mapping on a subset $C$ of a Banach space $X$. If $T$ satisfies the following property for all $x, y \in C$, then $T$ is said to satisfy condition ($C$) (($C$)–mapping).

$$\frac{1}{2}\|x - Tx\| \leq \|x - y\| \Rightarrow \|Tx - Ty\| \leq \|x - y\|$$

**Definition 3** ([1]) Let $T: C \to C$ be a mapping on a subset $C$ of a Banach space $X$. If there exists $a \in (1/2, 1)$ such that $T$ satisfies the following inequality for all $\alpha \in [a, 1], x \in C, y \in C(T, x, \alpha)$, then $T$ is said to satisfy condition ($D_a$) (($D_a$)–mapping).

$$\|Tx - Ty\| \leq \|x - y\|.$$

$$(C(T, x, \alpha) = \{y \in C \mid y = (1-\alpha)p + \alpha Tq, \ p, q \in C,$$
$$\|Tp - p\| \leq \|Tx - x\|, \|Tq - q\| \leq \|Tx - x\|\})$$

The ($D_a$)-mapping is weaker than nonexpansive mapping and stronger than quasinonexpansive mapping. And, it is different from the (C)-mapping. ([1])

The following lemma is useful to prove the existence and convergence theorems for ($D_a$)-mappings.

**Lemma 1** ([1]) Let $T: C \to C$ be a ($D_a$)-mapping. Then $T$ satisfies the following property for all $x, y \in C$.

$$\|Tx - x\| \leq \|Ty - y\| \Rightarrow \|x - Ty\| \leq 3\|Tx - x\| + \|x - y\|$$

**Definition 4** ([5]) Let $X$ be a Banach space. If the following inequality holds for any weakly convergent sequence $\{x_n\}$ and $y \neq x$ ($x$ is a weak limit), then $X$ is said to satisfy Opial's condition.
$$\liminf_{n\to\infty} \|x_n - x\| < \liminf_{n\to\infty} \|x_n - y\|$$

**Definition 5** ([5]) Let $\{x_n\}$ be a bounded sequence of a Banach space $X$. For $x \in C \subset X$, we set
$$r(x, \{x_n\}) := \limsup_{n\to\infty} \|x_n - x\|..$$

The asymptotic radius of $\{x_n\}$ relative to $C$ is defined by
$$r(C, \{x_n\}) := \inf \{r(x, \{x_n\}) | \ x \in C\}.$$

The asymptotic center of $\{x_n\}$ relative to $C$ is the set
$$A(C, \{x_n\}) := \{x \in C | \ r(x, \{x_n\}) = r(C, \{x_n\}) \}.$$

It is well known that $A(C, \{x_n\})$ consists only one point in uniformly convex Banach spaces.

**Lemma 2** ([7]) Let $X$ be a uniformly convex Banach space and for all $n \geq 1$.
$$0 < a \leq s_n \leq b < 1$$
Assume that two sequences $\{x_n\}, \{y_n\} \subset X$ satisfy the following conditions for some $d \geq 0$.
$$\limsup_{n\to\infty} \|x_n\| \leq d, \ \limsup_{n\to\infty} \|y_n\| \leq d, \ \limsup_{n\to\infty} \|s_n x_n + (1-s_n) y_n\| = d.$$

Then $\lim_{n\to\infty} \|x_n - y_n\| = 0$ holds.

## 2. Existence and convergence theorems

In this section we prove the existence for $(D_a)$–mappings using (*) and then establish the convergence theorems for (**).

Next lemmas are useful for the proof of the existence.

**Lemma 3** ([1]) Let $T$ be a $(D_a)$–mapping on a nonempty bounded convex subset $C$ of a Banach space $X$ and $\{x_n\}$ be a sequence generated by (*).

Then $\lim_{n\to\infty} \|x_n - Tx_n\| = 0$ holds.

**Lemma 4** Let $X$ be a Banach space with the Opial's property and $T$ be a $(D_a)$–mapping on a nonempty subset $C$ of $X$. If $\{x_n\} \subset C$ converges weakly to $z$ and $\lim_{n\to\infty} \|x_n - Tx_n\| = 0$ holds, then $Tz = z$.

**Proof.** Assume that $Tz \neq z$.

Since $\lim_{n\to\infty} \|x_n - Tx_n\| = 0$, by taking subsequences if necessary, we have
$$\|x_n - Tx_n\| \leq \|Tz - z\|$$
for all $n \in \mathbf{N}$.

Then from Lemma 1, we have

$$\|x_n - Tz\| \leq 3\|Tx_n - x_n\| + \|x_n - z\|.$$

Taking limits as $n \to \infty$,

$$\liminf_{n \to \infty} \|x_n - Tz\| \leq \liminf_{n \to \infty} \|x_n - z\|$$

holds.

On the other hand, from the Opial's property we have

$$\liminf_{n \to \infty} \|x_n - z\| < \liminf_{n \to \infty} \|x_n - Tz\|.$$

This is a contradiction, so $Tz = z$. □

**Theorem 1** Let $X$ be a Banach space with the Opial's property and $T$ be a ($D_a$)–mapping on a nonempty weakly compact convex subset $C$ of $X$. Then $T$ has a fixed point.

**Proof.** Since $C$ is weakly compact, it is bounded.

Let $\{x_n\}$ be a sequence generated by (*) and a subsequence $\{x_{n_k}\}$ converges weakly to $z \in C$.

Then from Lemma 3 and 4, we have $Tz = z$. This completes the proof. □

Now we establish the convergence theorems based on the existence of a fixed point.

Firstly, we prove some useful lemmas to prove the convergence.

**Lemma 5** Let $T$ be a ($D_a$)–mapping on a nonempty convex subset $C$ of a Banach space $X$. Assume that $\{x_n\}$ is a sequence generated by (**) and $T$ has a fixed point in $C$. Then for any fixed point $p$, $\lim_{n \to \infty} \|x_n - p\|$ exists.

**Proof.** The Proof is similar to the proof of Lemma 4 in [3]. □

**Lemma 6** Let $T$ be a ($D_a$)–mapping on a nonempty closed convex subset $C$ of a uniformly convex Banach space $X$. Assume that $\{x_n\}$ is a sequence generated by (**). Then $T$ has a fixed point if and only if $\{x_n\}$ is bounded and $\lim_{n \to \infty} \|x_n - Tx_n\| = 0$.

**Proof.** Proof of necessity is similar to the proof of Lemma 5 in [3].

Conversely, assume that $\{x_n\}$ is bounded and $\lim_{n \to \infty} \|x_n - Tx_n\| = 0$.

Let $p \in A(C, \{x_n\})$. Assume that $p \neq Tp$.

Then, similarly to the proof of Lemma 4, we have

$$\limsup_{n \to \infty} \|x_n - Tp\| \leq \limsup_{n \to \infty} \|x_n - p\|$$

Hence we obtain

$$r(Tp, \{x_n\}) = \limsup_{n \to \infty} \|x_n - Tp\| \leq \limsup_{n \to \infty} \|x_n - p\| = r(p, \{x_n\}) = r(C, \{x_n\})$$

, which leads to $Tp \in A(C, \{x_n\})$.

On the other hand, since $A(C, \{x_n\})$ is singleton in uniformly convex Banach spaces, we have $p = Tp$.

This is a contradiction, so $p$ is a fixed point. □

Now we establish weak and strong convergence theorems using (**).

**Theorem 2** Let $T$ be a ($D_a$)–mapping on a nonempty weakly compact convex subset $C$ of a uniformly convex Banach space $X$. Assume that $\{x_n\}$ is a sequence generated by (**). If $X$ satisfies Opial's condition, then $\{x_n\}$ converges weakly to a fixed point of $T$.

**Proof.** From Theorem 1, $T$ has a fixed point. Let $p \in F(T)$. From Lemma 4, $\lim_{n\to\infty} \|x_n - p\|$ exists.

Now we claim that any subsequence of $\{x_n\}$ has a unique weak limit. To prove this, let $x, y$ be weak limits of two subsequences $\{x_{n_j}\}, \{x_{n_k}\} \subset \{x_n\}$, respectively.

From Lemma 6, we have $\lim_{n\to\infty} \|x_n - Tx_n\| = 0$. So by Lemma 4, we have $x = Tx$, $y = Ty$.

Assume that $x \neq y$. Then from Opial's property, we obtain
$$\lim_{n\to\infty} \|x_n - x\| = \lim_{n_j \to\infty} \|x_{n_j} - x\|$$
$$< \lim_{n_j \to\infty} \|x_{n_j} - y\| = \lim_{n\to\infty} \|x_n - y\| = \lim_{n_k \to\infty} \|x_{n_k} - y\|.$$
$$< \lim_{n_k \to\infty} \|x_{n_k} - x\| = \lim_{n\to\infty} \|x_n - x\|$$

This is a contradiction, so we have $x = y$. □

**Theorem 3** Let $T$ be a ($D_a$)–mapping on a nonempty closed convex subset $C$ of a uniformly convex Banach space $X$ with the Opial's property. Assume that $\{x_n\}$ is a sequence generated by (**).

Then $\{x_n\}$ converges to a fixed point of $T$ if and only if
$$\liminf_{n\to\infty} d(x_n, F(T)) = 0 \text{ or } \limsup_{n\to\infty} d(x_n, F(T)) = 0$$

, where $d(x_n, F(T)) = \inf \{\|x_n - p\| : p \in F(T)\}$.

**Proof.** Proof is similar to the proof of Theorem 2 in [3]. □

**Theorem 4** Let $T$ be a ($D_a$)–mapping on a nonempty compact convex subset $C$ of a uniformly convex Banach space $X$. Assume that $\{x_n\}$ is a sequence generated by (**). Then $\{x_n\}$ converges to a fixed point of $T$.

**Proof.** From Theorem 3.7 in [1], $T$ has a fixed point. Then from Lemma 6, we have $\lim_{n\to\infty} \|x_n - Tx_n\| = 0$. Since $C$ is compact, there exists a subsequence $\{x_{n_k}\}$ of $\{x_n\}$ converging to $p \in C$.

Assume that $p$ is not a fixed point of $T$. Then, there exists a number $k_0$ such that

$$\|Tx_{n_k} - x_{n_k}\| \leq \|Tp - p\|$$

, for all $k \geq k_0$.

From Lemma 1, we have

$$\|x_{n_k} - Tp\| \leq 3\|Tx_{n_k} - x_{n_k}\| + \|x_{n_k} - p\|.$$

Since $\lim_{k \to \infty} \|x_{n_k} - p\| = 0$, $\lim_{n \to \infty} \|x_n - Tx_n\| = 0$, we have $\lim_{k \to \infty} \|x_{n_k} - Tp\| = 0$.

Therefore, we have $Tp = p$.

This is a contradiction, so $p$ is a fixed point of $T$. □

Now we introduce the following condition to prove the strong convergence theorem on weakly compact convex subsets.

**Definition 6** ([8]) A map $T: C \to C$ is said to satisfy condition ($I$), if there is a nondecreasing function $h: [0, \infty) \to [0, \infty)$ with $h(0) = 0$ and $h(r) > 0$ for any $r > 0$, such that
$$d(x, Tx) \geq h(d(x, F(T)))$$
, for all $x \in C$, where $d(x, F(T)) = \inf\{d(x, p): p \in F(T)\}$.

**Theorem 5** Let $T$ be a ($D_a$)–mapping on a nonempty weakly compact convex subset $C$ of a uniformly convex Banach space $X$ with the Opial's property. Assume that $\{x_n\}$ is a sequence generated by (**) and $T$ satisfies the condition ($I$). Then $\{x_n\}$ converges to a fixed point of $T$.

**Proof.** From Theorem 1, $T$ has a fixed point. And a weakly compact convex subset of a uniformly convex Banach space is a bounded closed convex subset.

The next process is similar to the proof of Theorem 4 in [3]. □

## 3. Numerical Example

In this section, we give an example to show that (**) converges to a fixed point faster than (*).

Define a mapping $T: [0, 1] \to \mathbf{R}$ by

$$T(x) = \begin{cases} \dfrac{x}{2}, & x \neq 1, \\ \dfrac{5}{8}, & x = 1. \end{cases}$$

Then, $T$ does not satisfy condition ($C$).

For $x = 1, y = 0.8$, we have

$$\frac{1}{2}\|x - Tx\| = \frac{3}{16} < \frac{1}{5} = \|x - y\|$$

but

$$\|Tx - Ty\| = \frac{9}{40} > \frac{1}{5} = \|x - y\|$$

, so it does not satisfy condition ($C$).

Now we verify that $T$ satisfies condition ($D_a$) for $a \geq \dfrac{2}{3}$.

Case 1: Let $x = 1$. Then for $y \in C(T, x, \alpha)$, we have $\dfrac{y}{2} \leq \dfrac{3}{8}$. So we obtain

$$\|Tx - Ty\| = \dfrac{5}{8} - \dfrac{y}{2} \leq 1 - y = \|x - y\|.$$

Case 2: Let $x \neq 1$.

If $y \in C(T, x, \alpha)$, $y \neq 1$, then we have $\dfrac{y}{2} \leq \dfrac{x}{2}$. So we have

$$\|Tx - Ty\| = \dfrac{x}{2} - \dfrac{y}{2} \leq x - y = \|x - y\|.$$

If $y \in C(T, x, \alpha)$, $y = 1$, then we have $\dfrac{3}{8} \leq \dfrac{x}{2}$. So we have

$$\|Tx - Ty\| = \dfrac{x}{2} - \dfrac{5}{8} \leq 0 \leq \|x - y\|.$$

Hence, $T$ satisfies condition ($D_a$).

With help of Mathematica Program Software, we obtain the comparison Table 1 and Figure 1 for (*) and (**) with $a_n = 0.85$, $b_n = 0.65$, $c_n = 0.45$, $x_1 = 0.9$.

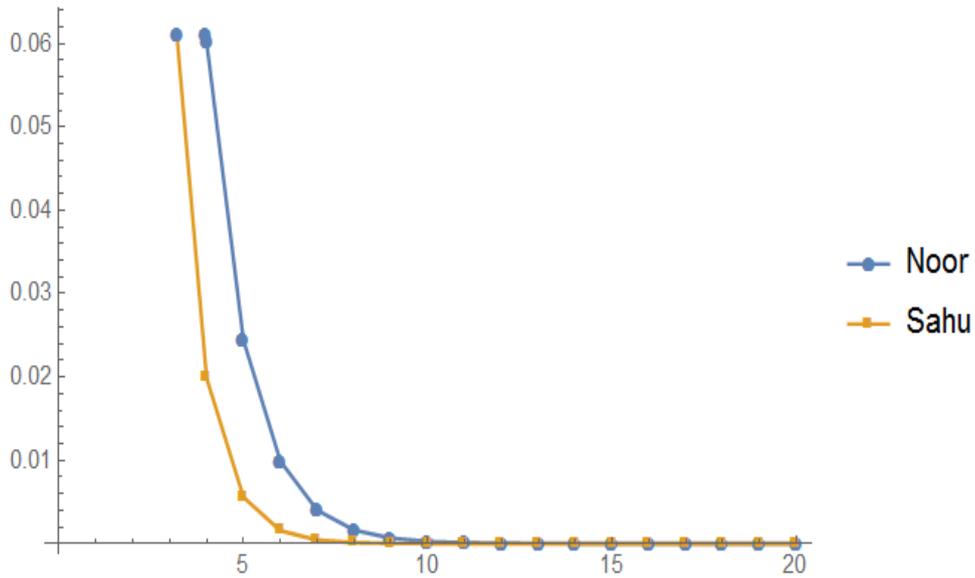

Figure 1

|    | Noor | Sahu |
|----|------|------|
| 1  | 0.9 | 0.9 |
| 2  | 0.365217 | 0.252408 |
| 3  | 0.148204 | 0.0707886 |
| 4  | 0.0601407 | 0.0198529 |
| 5  | 0.0244049 | 0.0055678 |
| 6  | 0.00990344 | 0.00156151 |
| 7  | 0.00401878 | 0.00043793 |
| 8  | 0.00163081 | 0.000122819 |
| 9  | 0.000661778 | 0.0000344449 |
| 10 | 0.000268547 | $9.66018 \times 10^{-6}$ |
| 11 | 0.000108976 | $2.70923 \times 10^{-6}$ |
| 12 | 0.000044222 | $7.59811 \times 10^{-7}$ |
| 13 | 0.0000179451 | $2.13091 \times 10^{-7}$ |
| 14 | $7.28208 \times 10^{-6}$ | $5.97621 \times 10^{-8}$ |
| 15 | $2.95505 \times 10^{-6}$ | $1.67605 \times 10^{-8}$ |
| 16 | $1.19915 \times 10^{-6}$ | $4.70053 \times 10^{-9}$ |
| 17 | $4.86611 \times 10^{-7}$ | $1.31828 \times 10^{-9}$ |
| 18 | $1.97465 \times 10^{-7}$ | $3.69715 \times 10^{-10}$ |
| 19 | $8.01307 \times 10^{-8}$ | $1.03688 \times 10^{-10}$ |
| 20 | $3.25168 \times 10^{-8}$ | $2.90796 \times 10^{-11}$ |

Table 1

From the figure and table, we can see that (**) converges faster than (*) to the fixed point $p = 0$.